\documentclass{article}
\usepackage{amscd}
\usepackage{amsmath}
\usepackage{amssymb}
\usepackage{amsthm}
\usepackage{ascmac}
\usepackage{graphicx}
\usepackage{here}
\usepackage{mathtools}\numberwithin{equation}{section}
\usepackage{minitoc}
\usepackage{multicol}
\usepackage{multirow}
\usepackage{setspace}
\theoremstyle{definition}

\newtheorem{q}{Problem}
\newtheorem{eg}{Example}
\newtheorem{rmk}{Remark}
\theoremstyle{theorem}
\newtheorem{thm}{Theorem}
\newtheorem{prop}{Proposition}
\newtheorem{lem}{Lemma}

\newenvironment{pf}{\begin{proof}}{\end{proof}}
\title{Diophantine equation related to angle bisectors\\
and solutions of Pell's equations}
\author{Takashi HIROTSU}
\date{\today}
\begin{document}
\maketitle
\begin{abstract}
It is important in drawing techniques to find combinations of two straight lines and their angle bisectors whose slopes are all rational numbers. 
This problem is reduced to solving the Diophantine equation \mbox{$(a-c)^2(b^2+1) = (b-c)^2(a^2+1).$} 
In this article, we describe all nontrivial integral solutions of the equation with solutions of negative Pell's equations. 
The formula is proven by certain properties of solutions of Pell's equations like those of half-companion Pell numbers and Pell numbers. 
We also give a formula for its rational solutions produced by Pythagorean triples with identical legs.
\end{abstract}
\section{Introduction}\label{sec-intro}
On the coordinate plane, the slope of the bisector of the acute angle between two straight lines with slopes $1$ and $7$ is $2.$ 
In addition, for example in the following cases, the slopes of the bisectors of the acute angles between two straight lines become integers.
\begin{table}[H]
\centering
\begin{tabular}{l|cccccccccccc}\hline 
\multirow{2}{*}{Slopes of two straight lines} & $1$ & $1$ & $2$ & $3$ & $4$ & $5$ & $6$ & $7$ & $7$ & $7$ & $8$ & $9$ \\ \cline{2-13} 
{} & $7$ & $-7$ &$38$ & $117$ & $268$ & $515$ & $882$ & $41$ & $1393$ & $-41$ & $2072$ & $2943$ \\ \hline 
The slope of an angle bisector & $2$ & $3$ & $4$ & $6$ & $8$ & $10$ & $12$ & $12$ & $14$ & $17$ & $16$ & $18$ \\ \hline 
\end{tabular}
\end{table}\par
In general, the following problem can be considered. 
Note that the two bisectors of the angles between given two straight lines are perpendicular to each other.
\begin{q}\label{q-bisec}
For which integers (or rational numbers) $a$ and $b$ is the slope of one of the angle bisectors between two straight lines with slopes $a$ and $b$ an integer (or a rational number)?
\end{q}
Essentially, Problem \ref{q-bisec} has the meaning when the bisectors of $\angle AOB$ can be drawn by connecting $O$ and other lattice points for given lattice points $O,$ $A,$ and $B,$ which is important in drawing techniques. 
By the following proposition, Problem \ref{q-bisec} is reduced to solving the equation 
\begin{equation} 
(a-c)^2(b^2+1) = (b-c)^2(a^2+1) \tag{$\star$} \label{star}  
\end{equation} 
in $\mathbb Z$ (or $\mathbb Q$).
\begin{prop}\label{prop-star}
Let $a,$ $b,$ $c \in \mathbb R$ with $|a| \neq |b|.$ 
If the slopes of the angle bisectors between two straight lines with slopes $a$ and $b$ are $c$ and $-c^{-1},$ then $a,$ $b,$ and $c$ satisfy \eqref{star}.
\end{prop}
\begin{pf}
It suffices to show the case when $y = cx$ is one of the angle bisectors between $y = ax$ and $y = bx.$  
Then the distances from a point $(t,ct) \neq (0,0)$ on $y = cx$ to $ax-y = 0$ and $bx-y = 0$ are equal to each other. 
This implies 
\[\frac{|at-ct|}{\sqrt{a^2+(-1)^2}} = \frac{|bt-ct|}{\sqrt{b^2+(-1)^2}},\] 
or equivalently, 
\[ |a-c|\sqrt{b^2+1} = |b-c|\sqrt{a^2+1}.\]  
Squaring both sides, we obtain \eqref{star}.
\end{pf}
\begin{rmk}\label{rmk-star}
Note that \eqref{star} is equivalent to $(ac+1)^2(b^2+1) = (bc+1)^2(a^2+1),$ which is derived from the same argument as above and the fact that the other angle bisector is $y = -c^{-1}x.$ 
Proposition \ref{prop-star} can also be proven by the addition formula of the tangent function, or the formula for the inner product of two vectors.
\end{rmk}\par
For solutions of \eqref{star}, the following properties are fundamental. 
We say that a solution $(a,b,c)$ of \eqref{star} is {\it trivial} if $|a| = |b|.$
\begin{prop}\label{prop-why}
\begin{enumerate}
\item[{\rm (1)}]
If $(a,b,c) = (a_1,b_1,c_1)$ is a real solution of \eqref{star}, then so are $(a,b,c) = (a_1,b_1,-c_1{}^{-1})$ $(c_1 \neq 0)$ and $(a,b,c) = (-a_1,-b_1,-c_1),$ $(b_1,a_1,c_1).$
\item[{\rm (2)}]
If $a = 0$ or $b = 0,$ then \eqref{star} has no nontrivial integral solutions.
\item[{\rm (3)}]
For any nontrivial integral solution $(a,b,c) = (a_1,b_1,c_1)$ of \eqref{star}, there exists a square-free integer $d > 1$ such that $a_1$ and $b_1$ are the $x$-components of integral solutions of $x^2-dy^2 = -1.$
\end{enumerate}
\end{prop}
\begin{pf}
\begin{enumerate}
\item[(1)]
See Remark \ref{rmk-star} for $(a,b,c) = (a_1,b_1,-c_1{}^{-1}).$ 
The others are obvious.
\item[(2)]
Let $a,$ $b \in \mathbb Z$ with $|a| \neq |b|.$ 
Solving \eqref{star} for $c,$ we obtain 
\[ c = \frac{ab-1\pm\sqrt{(a^2+1)(b^2+1)}}{a+b},\] 
which is not an integer if $a = 0$ or $b = 0,$ since $\sqrt{e^2+1} \notin \mathbb Q$ for any $e \in \mathbb Z\setminus\{ 0\}.$ 
This proves the desired assertion.
\item[(3)]
For any prime number $p,$ with respect to the normalized $p$-adic additive valuation $\mathrm{ord}_p:\mathbb Q^\times\to\mathbb Z,$ the parities of $\mathrm{ord}_p(a_1{}^2+1)$ and $\mathrm{ord}_p(b_1{}^2+1)$ coincide with each other, since $(a_1-c_1)^2$ and $(b_1-c_1)^2$ are square numbers. 
Let $d$ be the product of every prime number $p$ such that these valuations are odd, which is not equal to $1$ since $a_1{}^2+1$ and $b_1{}^2+1$ are not square numbers. 
Then there exist $a_2,$ $b_2 \in \mathbb Z$ such that $a_1{}^2+1 = da_2{}^2$ and $b_1{}^2+1 = db_2{}^2,$ which implies that $a_1$ and $b_1$ are the $x$-components of integral solutions of $x^2-dy^2 = -1.$\qedhere
\end{enumerate}
\end{pf}
\begin{rmk}
The {\it negative} Pell's equation $x^2-dy^2 = -1$ does not always have an integral solution, for example when $d = 3.$ 
It has integral solutions, for example when 
\[ d = 2,\ 5,\ 10,\ 13,\ 17,\ 26,\ 29,\ 37,\ 41,\ 53,\ 58,\ 61,\ 65,\ 73,\ 74,\ 82,\ 85,\ 89,\ 97.\] 
Such integers have no prime divisors congruent to $3$ modulo $4,$ but this condition is not sufficient. 
For example, $x^2-34y^2 = -1$ has no integral solutions.
\end{rmk}\par
The first main result of this article is the following formula for all nontrivial integral solutions of \eqref{star}, whose proof is given in Section \ref{sec-int} after preparation in Section \ref{sec-pell}. 
Let $d > 1$ be a square-free integer. 
For positive integral solutions $(x,y) = (a_1,a_2),$ $(b_1,b_2)$ of $|x^2-dy^2| = 1,$ we say that $(a_1,a_2)$ is {\it smaller} than $(b_1,b_2)$ if $a_2 < b_2.$
\begin{thm}\label{thm-int}
For each square-free integer $d > 1$ such that $x^2-dy^2 = -1$ has an integral solution, we denote the $n$-th smallest positive integral solution of $|x^2-dy^2| = 1$ by $(x,y) = (f_n^{(d)},g_n^{(d)}).$ 
Every nontrivial integral solution $(a,b,c)$ of \eqref{star} is given by 
\begin{align} 
(a,b,c) = &\pm\left( f_{(2m-1)(2n-1)}^{(d)},f_{(2m-1)(2n+1)}^{(d)},\frac{g_{(2m-1)\cdot 2n}^{(d)}}{g_{2m-1}^{(d)}}\right), \label{sol-d} \\
&\pm (f_{2n-1}^{(2)},-f_{2n+1}^{(2)},f_{2n}^{(2)}) \label{sol-2} 
\end{align} 
for some integers $d,$ $m,$ $n > 0,$ after switching $a$ and $b$ if necessary, where \eqref{sol-d} contains the case when $d = 2.$ 
Conversely, every triple $(a,b,c)$ of Form \eqref{sol-d} or \eqref{sol-2} is an integral solution of \eqref{star}. 
\end{thm}
\begin{eg}
For a given integer $e > 0,$ \eqref{star} has integral solutions 
\[ (a,b,c) = \pm (e,e(4e^2+3),2e),\] 
where $e$ and $e(4e^2+3)$ are the $x$-components of the first and third smallest solutions $(x,y) = (e,1)$ and $(x,y) = (e(4e^2+3),4e^2+1)$ of $|x^2-(e^2+1)y^2| = 1,$ respectively.
\end{eg}
\begin{eg}
Integral solutions of $|x^2-2y^2| = 1$ produce the following integral solutions of \eqref{star}: 
the first few solutions of Form \eqref{sol-d} are
\begin{table}[H]
\centering
\begin{tabular}{lll} 
$(f_1^{(2)},f_3^{(2)},g_2^{(2)}/g_1^{(2)})$ & $(f_3^{(2)},f_5^{(2)},g_4^{(2)}/g_1^{(2)})$ & $(f_5^{(2)},f_7^{(2)},g_6^{(2)}/g_1^{(2)})$ \\ 
$= (1,7,2),$ & $= (7,41,12),$ & $= (41,239,70),$ \\[2mm] 
$(f_3^{(2)},f_9^{(2)},g_6^{(2)}/g_3^{(2)})$ & $(f_9^{(2)},f_{15}^{(2)},g_{12}^{(2)}/g_3^{(2)})$ & $(f_{15}^{(2)},f_{21}^{(2)},g_{18}^{(2)}/g_3^{(2)})$ \\ 
$= (7,1393,14),$ & $= (1393,275807,2772),$ & $= (275807,54608393,548842),$ \\[2mm] 
$(f_5^{(2)},f_{15}^{(2)},g_{10}^{(2)}/g_5^{(2)})$ & $(f_{15}^{(2)},f_{25}^{(2)},g_{20}^{(2)}/g_5^{(2)})$ & $(f_{25}^{(2)},f_{35}^{(2)},g_{30}^{(2)}/g_5^{(2)})$ \\ 
$= (41,275807,82),$ & $= (275807,1855077841,551532),$ & $= (1855077841,12477253282759,3709604150),$ 
\end{tabular}
\end{table}
and the first few solutions of Form \eqref{sol-2} are
\begin{table}[H]
\centering
\begin{tabular}{lllll} 
$(f_1^{(2)},-f_3^{(2)},f_2^{(2)})$ & $(f_3^{(2)},-f_5^{(2)},f_4^{(2)})$ & $(f_5^{(2)},-f_7^{(2)},f_6^{(2)})$ & $(f_7^{(2)},-f_9^{(2)},f_8^{(2)})$ & $(f_9^{(2)},-f_{11}^{(2)},f_{10}^{(2)})$ \\ 
$= (1,-7,3),$ & $= (7,-41,17),$ & $= (41,-239,99),$ & $= (239,-1393,577),$ & $= (1393,-8119,3363).$ 
\end{tabular}
\end{table}
\end{eg}
\begin{eg}
Integral solutions of $|x^2-5y^2| = 1$ produce the following integral solutions of \eqref{star}: 
the first few solutions of Form \eqref{sol-d} are
\begin{table}[H]
\centering
\begin{tabular}{lll} 
$(f_1^{(5)},f_3^{(5)},g_2^{(5)}/g_1^{(5)})$ & $(f_3^{(5)},f_5^{(5)},g_4^{(5)}/g_1^{(5)})$ & $(f_5^{(5)},f_7^{(5)},g_6^{(5)}/g_1^{(5)})$ \\ 
$= (2,38,4),$ & $= (38,682,72),$ & $= (682,12238,1292),$ \\[2mm] 
$(f_3^{(5)},f_9^{(5)},g_6^{(5)}/g_3^{(5)})$ & $(f_9^{(5)},f_{15}^{(5)},g_{12}^{(5)}/g_3^{(5)})$ &  \\ 
$= (38,219602,76),$ & $= (219602,1268860318,439128),$ &  \\[2mm] 
$(f_5^{(5)},f_{15}^{(5)},g_{10}^{(5)}/g_5^{(5)})$  &  \\ 
$= (682,1268860318,1364).$ &  & 
\end{tabular}
\end{table}
\end{eg}
The second main result is the following formula for a kind of rational solutions of \eqref{star}, whose proof is given in Section \ref{sec-rat}.
\begin{thm}\label{thm-rat}
Let $w$ be a multiple of $4$ greater than $4,$ the twice of an odd composite number, or an odd composite number. 
Then $u^2-v^2 = -w^2$ has distinct positive integral solutions $(u,v) = (x_1,x_2),$ $(y_1,y_2),$ and \eqref{star} has rational solutions of the form 
\begin{equation} 
(a,b,c) = \pm\left(\frac{x_1}{w},\frac{y_1}{w},\frac{x_1y_2+x_2y_1}{w(y_2+x_2)}\right),\ \pm\left(\frac{x_1}{w},\frac{y_1}{w},\frac{x_1y_2-x_2y_1}{w(y_2-x_2)}\right). \label{sol-rat} 
\end{equation} 
\end{thm}
\begin{eg}
The pair of Pythagorean triples $(5,12,13)$ and $(35,12,37)$ with identical leg $12$ produces the rational solutions 
\[\left(\frac{5}{12},\frac{35}{12},\frac{5\cdot 37+13\cdot 35}{12(37+13)}\right) = \left(\frac{5}{12},\frac{35}{12},\frac{16}{15}\right), \quad \left(\frac{5}{12},\frac{35}{12},\frac{5\cdot 37-13\cdot 35}{12(37-13)}\right) = \left(\dfrac{5}{12},\dfrac{35}{12},-\dfrac{15}{16}\right)\] 
of \eqref{star}.
\end{eg}
\section{Properties of Solutions of Pell's Equations}\label{sec-pell}
Let $d > 1$ be a square-free integer. 
Assume that $x^2-dy^2 = -1$ has an integral solution. 
In the real quadratic field $\mathbb Q(\sqrt d),$ we denote the conjugate $a_1-a_2\sqrt d$ of $\alpha = a_1+a_2\sqrt d$ over $\mathbb Q$ by $\alpha ',$ where $a_1,$ $a_2 \in \mathbb Q.$ 
In this section, we denote $f_n^{(d)}$ and $g_n^{(d)}$ defined in Theorem \ref{thm-int} by $f_n$ and $g_n,$ respectively, without the indices $(d).$ 
For convenience, let $f_0 = 1$ and $g_0 = 0.$ 
In the case when $d = 2,$ the terms of $(f_n)$ and $(g_n)$ are known as {\itshape half-companion Pell numbers}, {\itshape Pell numbers}, respectively. 
Furthermore, let $\varepsilon = f_1+g_1\sqrt d.$ 
Note that $\varepsilon\varepsilon ' = -1$ by assumption.\par
In this section, we describe the properties of $(f_n)$ and $(g_n)$ used in the proof of Theorem \ref{thm-int} and certain related properties. 
Since most of them are already known (see \cite{car}), we give brief proofs of them. 
The following two propositions are well-known (see \cite[Section 2.4]{kos}).
\begin{prop}\label{prop-gen-pel}
The sequences $(f_n)$ and $(g_n)$ are strictly increasing, and satisfy 
\begin{equation} 
\varepsilon ^n = f_n+g_n\sqrt d \label{unit-sum} 
\end{equation} 
and 
\begin{equation} 
f_n{}^2-dg_n{}^2 = (-1)^n. \label{pell-eq} 
\end{equation} 
Their general terms are given by 
\begin{align} 
f_n &= \frac{\varepsilon ^n+\varepsilon '^n}{2}, \label{seq-f} \\ 
g_n &= \frac{\varepsilon ^n-\varepsilon '^n}{2\sqrt d}. \label{seq-g} 
\end{align} 
\end{prop}
\begin{pf}
The increasingness follows from the definition.\par
For each integer $n > 0,$ let $\varepsilon ^n = u_n+v_n\sqrt d$ with $u_n,$ $v_n \in \mathbb Z.$ 
Let $(x,y)$ be a positive integral solution of $|x^2-dy^2| = 1.$ 
Since $x+y\sqrt d \geq \varepsilon > 1,$ there exists an integer $n > 0$ such that 
\[\varepsilon ^n \leq x+y\sqrt d < \varepsilon ^{n+1},\] 
or equivalently, 
\[ 1 \leq (x+y\sqrt d)\varepsilon ^{-n} < \varepsilon.\] 
Since $\varepsilon ^{-n} = (-\varepsilon ')^n = (-1)^n(\varepsilon ^n)' = (-1)^n(u_n-v_n\sqrt d),$ the integers $x_0 = (-1)^n(xu_n-dyv_n)$ and $y_0 = (-1)^n(-xv_n+yu_n)$ satisfy $x_0+y_0\sqrt d = (x+y\sqrt d)\varepsilon ^{-n},$ $x_0-y_0\sqrt d = (x-y\sqrt d)(-\varepsilon )^n,$ and therefore 
\[ x_0{}^2-dy_0{}^2 = x^2-dy^2 = \pm 1.\] 
This implies $x_0+y_0\sqrt d = 1,$ $\varepsilon ^n = x+y\sqrt d,$ and therefore \eqref{unit-sum} by the increasingness of $(f_n)$ and $(g_n).$\par
Taking the conjugates in \eqref{unit-sum}, we obtain 
\begin{equation} 
\varepsilon '^n = f_n-g_n\sqrt d. \tag*{(2.1)${}'$} \label{unit-diff} 
\end{equation} 
Multiplying \eqref{unit-sum} and \ref{unit-diff}, we obtain \eqref{pell-eq}. 
Solving \eqref{unit-sum} and \ref{unit-diff} for $f_n$ and $g_n,$ we obtain \eqref{seq-f} and \eqref{seq-g}.
\end{pf}
\begin{prop}\label{prop-mag}
We have $f_n \geq g_n,$ where the equality holds if and only if $d = 2$ and $n = 1.$ 
We also have $f_1 \geq \sqrt{d-1}.$
\end{prop}
\begin{pf}
By \eqref{pell-eq}, we have 
\begin{align*} 
f_{2n-1}{}^2 &= dg_{2n-1}{}^2-1 \geq g_{2n-1}{}^2 \quad (\because (d-1)g_{2n-1}{}^2 \geq 1), \\ 
f_{2n}{}^2 &= dg_{2n}{}^2+1 > g_{2n}{}^2, 
\end{align*} 
and therefore $f_n \geq g_n.$ 
Furthermore, we have 
\begin{equation} 
f_1 = \sqrt{dg_1{}^2-1} \geq \sqrt{d-1}. \tag*{\qedhere} 
\end{equation} 
\end{pf}
The following proposition is also well-known and generalized for the Lucas sequences (see \cite[Chapter 2, IV]{rib}).
\begin{prop}
Let $m,$ $n \in \mathbb Z$ with $m \geq n \geq 0.$
\begin{enumerate}
\item[{\rm (1)}]
Addition formulas: we have 
\begin{align} 
f_{m+n} &= f_mf_n+dg_mg_n, \label{sum-f} \\ 
g_{m+n} &= f_mg_n+g_mf_n, \label{sum-g} \\ 
f_{m-n} &= (-1)^n(f_mf_n-dg_mg_n), \label{diff-f} \\ 
g_{m-n} &= (-1)^{n+1}(f_mg_n-g_mf_n). \label{diff-g} 
\end{align} 
\item[{\rm (2)}]
Double formulas: we have 
\begin{align} 
f_{2n} &= f_n^2+dg_n^2, \label{dbl-f} \\ 
g_{2n} &= 2f_ng_n. \label{dbl-g} 
\end{align} 
\end{enumerate}
\end{prop}
\begin{pf}
\begin{enumerate}
\item[(1)]
Describing $\varepsilon ^{m+n} = \varepsilon ^m\varepsilon ^n$ and $\varepsilon ^{m-n} = (-1)^n\varepsilon ^m\varepsilon '^n$ with terms of $(f_n)$ and $(g_n)$ by \eqref{unit-sum}, we obtain 
\begin{align*} 
f_{m+n}+g_{m+n}\sqrt d &= (f_m+g_m\sqrt d)(f_n+g_n\sqrt d) \\ 
&= (f_mf_n+dg_mg_n)+(f_mg_n+g_mf_n)\sqrt d, \\ 
f_{m-n}+g_{m-n}\sqrt d &= (-1)^n(f_m+g_m\sqrt d)(f_n-g_n\sqrt d) \\ 
&= (-1)^n(f_mf_n-dg_mg_n)+(-1)^{n+1}(f_mg_n-g_mf_n)\sqrt d. 
\end{align*} 
Comparing both sides, we obtain the desired identities, since $1$ and $\sqrt d$ are linearly independent over $\mathbb Q.$
\item[(2)]
Letting $m = n$ in \eqref{sum-f} and \eqref{sum-g}, we obtain the desired identities.\qedhere
\end{enumerate}
\end{pf}\par
The divisibility in $(g_n)$ depends only on that of the indices (see \cite[Theorem IV]{car} and \cite[Theorem 8.4]{kos}). 
A certain divisibility property in $(f_n)$ can be proven in a similar way (see \cite[Theorem V]{car}). 
These are summarized as follows. 
We denote the greatest common divisor of $a,$ $b \in \mathbb Z\setminus\{ 0\}$ by $\mathrm{gcd}\,(a,b).$
\begin{thm}\label{thm-div}
Let $m,$ $n \in \mathbb Z$ with $m \geq n > 0.$
\begin{enumerate}
\item[{\rm (1)}]
We have $\mathrm{gcd}\,(d,f_n) = \mathrm{gcd}\,(f_n,g_n) = 1.$
\item[{\rm (2)}]
If $m$ is a multiple of $n$ whose quotient is even, then $\mathrm{gcd}\,(f_m,f_n) = 1.$
\item[{\rm (3)}]
The following conditions are equivalent.
\begin{enumerate}
\item[{\rm (f1)}]
$f_m$ is a multiple of $f_n.$
\item[{\rm (f2)}]
Either we have $d = 2$ and $n = 1,$ or $m$ is a multiple of $n$ whose quotient is odd.
\end{enumerate}
\item[{\rm (4)}]
The following conditions are equivalent.
\begin{enumerate}
\item[{\rm (g1)}]
$g_m$ is a multiple of $g_n.$
\item[{\rm (g2)}]
$m$ is a multiple of $n.$
\end{enumerate}
\end{enumerate}
\end{thm}
\begin{pf}
We give the proof in the following order: (1), (f2) $\Rightarrow$ (f1) of (3), (2), (f1) $\Rightarrow$ (f2) of (3), (4).
\begin{enumerate}
\item[(1)]
Identity \eqref{pell-eq} implies $\mathrm{gcd}\,(d,f_n) = \mathrm{gcd}\,(f_n,g_n) = 1.$
\item[(3)]
We prove (f2) $\Rightarrow$ (f1). 
\begin{itemize}
\item
{\it Case 1}: Suppose that $d = 2$ and $n = 1.$ 
Then $f_m$ is a multiple of $f_n = 1,$ since $\varepsilon = 1+\sqrt 2.$
\item
{\it Case 2}: Suppose that $m = nq$ with an odd integer $q > 0.$ 
Since 
\[\frac{\varepsilon ^m+\varepsilon '^m}{2} = \frac{\varepsilon ^n+\varepsilon '^n}{2}\sum_{i = 0}^{q-1}(-1)^i\varepsilon ^{n(q-1-i)}\varepsilon '^{ni},\] 
we have 
\[\frac{f_m}{f_n} = \sum_{i = 0}^{q-1}(-1)^i\varepsilon ^{n(q-1-i)}\varepsilon '^{ni} \in \mathbb Z[\sqrt d]\cap\mathbb Q = \mathbb Z\] 
by \eqref{seq-f}. 
This implies that $f_m$ is a multiple of $f_n.$
\end{itemize}
\item[(2)]
Suppose that $m = nq$ with an even integer $q > 0.$ 
By \eqref{sum-f}, we have 
\[ f_m = f_{n(q-1)+n} = f_{n(q-1)}f_n+dg_{n(q-1)}g_n.\] 
Since $d$ and $g_n$ are coprime with $f_n$ by (1), we have $\mathrm{gcd}\,(f_m,f_n) = \mathrm{gcd}\,(g_{n(q-1)},f_n),$ and therefore this is a common divisor of $g_{n(q-1)}$ and $f_{n(q-1)}$ by (f2) $\Rightarrow$ (f1), which implies $\mathrm{gcd}\,(f_m,f_n) = 1$ by (1).
\item[(3)]
We prove (f1) $\Rightarrow$ (f2). 
Suppose that $f_m$ is a multiple of $f_n,$ and $d \neq 2$ or $n \neq 1.$ 
By Proposition \ref{prop-mag}, we have $f_n > 1.$ 
Let $q$ and $r$ be the quotient and remainder, respectively, when dividing $m$ by $n.$ 
By \eqref{sum-f}, we have 
\[ f_m = f_{nq+r} = f_{nq}f_r+dg_{nq}g_r.\]\par
Assume that $q$ is even, and let $q = 2^ek,$ where $e,$ $k > 0$ are integers and $k$ is odd. 
By \eqref{dbl-g}, we have $g_{nq} = 2f_{nq/2}g_{nq/2}.$ 
Repeating this $e$ times, we see that $g_{nq}$ is a multiple of $f_{nk}.$ 
By (f2) $\Rightarrow$ (f1), $f_{nk}$ is a multiple of $f_n.$ 
These imply that $g_{nq}$ is a multiple of $f_n.$ 
Since $f_{nq}$ is coprime with $f_n$ by (2), and $f_r$ is not a multiple of $f_n$ by $0 < f_r < f_n,$ we see that $f_m$ is not a multiple of $f_n.$
This is a contradiction.\par
Therefore $q$ is odd. 
Since $f_{nq}$ is a multiple of $f_n$ by (f2) $\Rightarrow$ (f1), $dg_{nq}g_r$ is a multiple of $f_n.$ 
Furthermore, $\mathrm{gcd}\,(f_n,g_{nq})$ is a common divisor of $f_{nq}$ and $g_{nq},$ which implies $\mathrm{gcd}\,(f_n,g_{nq}) = 1$ by (1). 
Since $d$ and $g_{nq}$ are coprime with $f_n,$ we see that $g_r$ is a multiple of $f_n.$ 
Since $0 \leq g_r < g_n < f_n$ by Propositions \ref{prop-gen-pel} and \ref{prop-mag}, we have $r = 0$ and $m = nq,$ which implies that $m$ is a multiple of $n$ whose quotient is odd.
\item[(4)]
We prove (g2) $\Rightarrow$ (g1). 
Suppose that $m = nq$ with an integer $q > 0.$ 
Since 
\[\frac{\varepsilon ^m-\varepsilon '^m}{2\sqrt d} = \frac{\varepsilon ^n-\varepsilon '^n}{2\sqrt d}\sum_{i = 0}^{q-1}\varepsilon ^{n(q-1-i)}\varepsilon '^{ni},\] 
we have 
\[\frac{g_m}{g_n} = \sum_{i = 0}^{q-1}\varepsilon ^{n(q-1-i)}\varepsilon '^{ni} \in \mathbb Z[\sqrt d]\cap\mathbb Q = \mathbb Z\] 
by \eqref{seq-g}. 
This implies that $g_m$ is a multiple of $g_n.$\par
We prove (g1) $\Rightarrow$ (g2). 
Suppose that $g_m$ is a multiple of $g_n.$ 
Let $q$ and $r$ be the quotient and remainder, respectively, when dividing $m$ by $n.$ 
By \eqref{sum-g}, we have 
\[ g_m = g_{nq+r} = f_{nq}g_r+g_{nq}f_r.\] 
Since $g_{nq}$ is a multiple of $g_n$ by (g2) $\Rightarrow$ (g1), $f_{nq}g_r$ is a multiple of $g_n.$ 
Furthermore, $\mathrm{gcd}\,(f_{nq},g_n)$ is a common divisor of $f_{nq}$ and $g_{nq},$ which implies $\mathrm{gcd}\,(f_{nq},g_n) = 1$ by (1). 
Since $f_{nq}$ is coprime with $g_n,$ we see that $g_r$ is a multiple of $g_n.$ 
Since $0 \leq g_r < g_n,$ we have $r = 0$ and $m = nq,$ which implies that $m$ is a multiple of $n.$\qedhere
\end{enumerate}
\end{pf}\par
The following formulas enable us to convert sums into products in $(f_n)$ and $(g_n).$
\begin{prop}
Let $m,$ $n \in \mathbb Z$ with $m > n \geq 0.$ 
Then we have 
\allowdisplaybreaks[1]
\begin{align} 
f_{m+n}+f_{m-n} &= \begin{cases} 
2f_mf_n & \text{if}\quad n \equiv 0\ (\mathrm{mod}\ 2), \\ 
2dg_mg_n & \text{if}\quad n \equiv 1\ (\mathrm{mod}\ 2), 
\end{cases} \label{sf} \\ 
f_{m+n}-f_{m-n} &= \begin{cases} 
2dg_mg_n & \text{if}\quad n \equiv 0\ (\mathrm{mod}\ 2), \\ 
2f_mf_n & \text{if}\quad n \equiv 1\ (\mathrm{mod}\ 2), 
\end{cases} \label{df} \\ 
g_{m+n}+g_{m-n} &= \begin{cases} 
2g_mf_n\;\, & \text{if}\quad n \equiv 0\ (\mathrm{mod}\ 2), \\ 
2f_mg_n & \text{if}\quad n \equiv 1\ (\mathrm{mod}\ 2), 
\end{cases} \label{sg} \\ 
g_{m+n}-g_{m-n} &= \begin{cases} 
2f_mg_n\;\, & \text{if}\quad n \equiv 0\ (\mathrm{mod}\ 2), \\ 
2g_mf_n & \text{if}\quad n \equiv 1\ (\mathrm{mod}\ 2). 
\end{cases} \label{dg} 
\end{align} 
\end{prop}
\begin{pf}
For any $n \in \mathbb Z,$ we have $\varepsilon ^n\varepsilon '^n = (\varepsilon\varepsilon ')^n = (-1)^n.$
By \eqref{seq-f} and \eqref{seq-g}, we have 
\begin{align*} 
2f_mf_n &= \frac{(\varepsilon ^m+\varepsilon '^m)(\varepsilon ^n+\varepsilon '^n)}{2} = \frac{(\varepsilon ^{m+n}+\varepsilon '^{m+n})+(-1)^n(\varepsilon ^{m-n}+\varepsilon '^{m-n})}{2} = f_{m+n}+(-1)^nf_{m-n}, \\ 
2dg_mg_n &= \frac{(\varepsilon ^m-\varepsilon '^m)(\varepsilon ^n-\varepsilon '^n)}{2} = \frac{(\varepsilon ^{m+n}+\varepsilon '^{m+n})+(-1)^{n+1}(\varepsilon ^{m-n}+\varepsilon '^{m-n})}{2} = f_{m+n}+(-1)^{n+1}f_{m-n}, \\ 
2f_mg_n &= \frac{(\varepsilon ^m+\varepsilon '^m)(\varepsilon ^n-\varepsilon '^n)}{2\sqrt d} = \frac{(\varepsilon ^{m+n}-\varepsilon '^{m+n})+(-1)^{n+1}(\varepsilon ^{m-n}-\varepsilon '^{m-n})}{2\sqrt d} = g_{m+n}+(-1)^{n+1}g_{m-n}, \\ 
2g_mf_n &= \frac{(\varepsilon ^m-\varepsilon '^m)(\varepsilon ^n+\varepsilon '^n)}{2\sqrt d} = \frac{(\varepsilon ^{m+n}-\varepsilon '^{m+n})+(-1)^n(\varepsilon ^{m-n}-\varepsilon '^{m-n})}{2\sqrt d} = g_{m+n}+(-1)^ng_{m-n}, 
\end{align*} 
which imply the desired identities.
\end{pf}
\section{Integral Solutions}\label{sec-int}
In this section, we prove Theorem \ref{thm-int}. 
\begin{lem}\label{lem-can}
Let $d > 1$ be a square-free integer. 
If $x^2-dy^2 = -1$ has integral solutions $(x,y) = (a_1,a_2),$ $(b_1,b_2)$ such that $|a_1| \neq |b_1|,$ then \eqref{star} has the rational solutions 
\begin{align*} 
(a,b,c) = \pm\left( a_1,b_1,\frac{a_1b_2+a_2b_1}{b_2+a_2}\right),\ \pm\left( a_1,b_1,\frac{a_1b_2-a_2b_1}{b_2-a_2}\right). 
\end{align*} 
\end{lem}
\begin{pf}
Let $a = a_1,$ $b = b_1,$ and $a_1{}^2-da_2{}^2 = b_1{}^2-db_2{}^2 = -1$ with $|a_1| \neq |b_1|.$ 
Then \eqref{star} implies 
\[ (a_1-c)^2db_2{}^2 = (b_1-c)^2da_2{}^2\] 
and therefore 
\[ (a_1-c)^2b_2{}^2 = (b_1-c)^2a_2{}^2.\] 
Solving for $c,$ we obtain 
\begin{equation} 
c = \frac{a_1b_2+a_2b_1}{b_2+a_2},\ \frac{a_1b_2-a_2b_1}{b_2-a_2}, \tag*{\qedhere} 
\end{equation} 
since $b_2{}^2-a_2{}^2 = (b_1{}^2-a_1{}^2)/d \neq 0.$
\end{pf}
\begin{pf}[{\bfseries Proof of Theorem \ref{thm-int}}]
Every nontrivial integral solution $(a,b,c)$ of \eqref{star} can necessarily be written in the form 
\[ (a,b,c) = \left( a_1,b_1,\frac{a_1b_2+a_2b_1}{b_2+a_2}\right),\ \left( a_1,b_1,\frac{a_1b_2-a_2b_1}{b_2-a_2}\right)\] 
for some $a_1,$ $a_2,$ $b_1,$ $b_2 \in \mathbb Z$ such that $a_1{}^2-da_2{}^2 = b_1{}^2-db_2{}^2 = -1$ by Proposition \ref{prop-why} and Lemma \ref{lem-can}. 
Henceforth, we consider the condition that the rational numbers 
\[ c_+ = \frac{a_1b_2+a_2b_1}{b_2+a_2} \quad\text{or}\quad c_- = \frac{a_1b_2-a_2b_1}{b_2-a_2}\] 
are integers. 
In the general case, we denote $f_n^{(d)}$ and $g_n^{(d)}$ by $f_n$ and $g_n,$ respectively. 
Suppose that $0 < a_1$ and $0 < a_2 < b_2.$ 
Then $(a_1,a_2) = (f_k,g_k)$ and $(b_1,b_2) = (\pm f_l,g_l)$ for some odd indices $k$ and $l$ such that $k < l.$
\begin{itemize}
\item
{\it Case 1}: Suppose that 
\[ (a_1,a_2) = (f_{2i-1},g_{2i-1}) \quad\text{and}\quad (b_1,b_2) = (f_{(2i-1)+(4j-2)},g_{(2i-1)+(4j-2)})\] 
for some integers $i,$ $j > 0.$ 
Then we have 
\[ c_+ = \frac{g_{4(i+j-1)}}{2f_{2(i+j-1)}g_{2j-1}} = \frac{g_{2(i+j-1)}}{g_{2j-1}}\] 
by \eqref{sum-g}, \eqref{sg}, and \eqref{dbl-g}. 
This implies that $c_+ \in \mathbb Z$ holds if and only if \mbox{$g_{2j-1} \mid g_{2(i+j-1)},$} which is equivalent to \mbox{$2j-1\mid 2(i+j-1)$} by Theorem \ref{thm-div} (4), and to \mbox{$2j-1 \mid 2i-1.$} 
Letting $j = m$ and $2i-1 = (2m-1)(2n-1),$ we obtain 
\[ (a_1,b_1,c_+) = \left( f_{(2m-1)(2n-1)},f_{(2m-1)(2n+1)},\frac{g_{(2m-1)\cdot 2n}}{g_{2m-1}}\right).\] 
Furthermore, $c_- = -c_+{}^{-1} \notin \mathbb Z,$ since $g_{(2m-1)\cdot 2n} > g_{2m-1}.$
\item
{\it Case 2}: Suppose that 
\[ (a_1,a_2) = (f_{2i-1},g_{2i-1}) \quad\text{and}\quad (b_1,b_2) = (-f_{(2i-1)+(4j-2)},g_{(2i-1)+(4j-2)})\] 
for some integers $i,$ $j > 0.$ 
Then we have 
\[ c_- = \frac{g_{4(i+j-1)}}{2g_{2(i+j-1)}f_{2j-1}} = \frac{f_{2(i+j-1)}}{f_{2j-1}}\] 
by \eqref{sum-g}, \eqref{dg}, and \eqref{dbl-g}. 
This implies that $c_- \in \mathbb Z$ holds if and only if \mbox{$f_{2j-1} \mid f_{2(i+j-1)},$} which is equivalent to $d = 2$ and $j = 1$ by Theorem \ref{thm-div} (3), since $2(i+j-1)$ is not a multiple of $2j-1$ whose quotient is odd. 
Letting $i = n,$ we obtain 
\[ (a_1,b_1,c_-) = ( f_{2n-1}^{(2)},-f_{2n+1}^{(2)},f_{2n}^{(2)}).\] 
Furthermore, $c_+ = -c_-{}^{-1} \notin \mathbb Z,$ since $f_{2n}^{(2)} > 1.$
\item
{\it Case 3}: Suppose that 
\[ (a_1,a_2) = (f_{2i-1},g_{2i-1}) \quad\text{and}\quad (b_1,b_2) = (f_{(2i-1)+4j},g_{(2i-1)+4j})\] 
for some integers $i,$ $j > 0.$ 
Then we have 
\[ c_+ = \frac{g_{4(i+j)-2}}{2g_{2(i+j)-1}f_{2j}} = \frac{f_{2(i+j)-1}}{f_{2j}}\] 
by \eqref{sum-g}, \eqref{sg}, and \eqref{dbl-g}. 
In addition, $f_{2(i+j)-1}$ is not a multiple of $f_{2j}$ by Theorem \ref{thm-div} (3), since $2(i+j)-1$ is not a multiple of $2j$ whose quotient is odd. 
These imply $c_+ \notin \mathbb Z.$ 
Furthermore, $c_- = -c_+{}^{-1} \notin \mathbb Z,$ since $f_{2(i+j)-1} > f_{2j}.$
\item
{\it Case 4}: Suppose that 
\[ (a_1,a_2) = (f_{2i-1},g_{2i-1}) \quad\text{and}\quad (b_1,b_2) = (-f_{(2i-1)+4j},g_{(2i-1)+4j})\] 
for some integers $i,$ $j > 0.$ 
Then we have 
\[ c_- = \frac{g_{4(i+j)-2}}{2f_{2(i+j)-1}g_{2j}} = \frac{g_{2(i+j)-1}}{g_{2j}}\] 
by \eqref{sum-g}, \eqref{dg}, and \eqref{dbl-g}. 
In addition, $g_{2(i+j)-1}$ is not a multiple of $g_{2j}$ by Theorem \ref{thm-div} (4), since $2(i+j)-1$ is not a multiple of $2j.$ 
These imply $c_- \notin \mathbb Z.$ 
Furthermore, $c_+ = -c_-{}^{-1} \notin \mathbb Z,$ since $g_{2(i+j)-1} > g_{2j}.$
\end{itemize}
Considering the sign changes, it is concluded that every nontrivial integral solution of \eqref{star} is given by \eqref{sol-d} or \eqref{sol-2} after switching $a$ and $b$ if necessary.\par
It is easy to verify that every triple of Form \eqref{sol-d} or \eqref{sol-2} is an integral solution of \eqref{star}.
\end{pf}
\section{Rational Solutions}\label{sec-rat}
In this section, we prove Theorem \ref{thm-rat} by using the following known result.
\begin{lem}[{\cite[Theorem 6]{tri}}]\label{lem-pyt}
Let $w > 0$ be an integer whose factorization into prime numbers is 
\[ w = 2^{e_0}p_1{}^{e_1}\cdots p_r{}^{e_r},\] 
where $p_1,$ $\ldots,$ $p_r$ are distinct odd prime numbers and $e_0,$ $e_1,$ $\ldots,$ $e_r \geq 0$ are integers. 
Then the number of positive integral solutions of $u^2-v^2 = -w^2$ is given by 
\begin{equation} 
\begin{cases} 
\dfrac{(2e_0-1)(2e_1+1)\cdots (2e_r+1)-1}{2} & \text{if }w\text{ is even}, \\ 
\dfrac{(2e_1+1)\cdots (2e_r+1)-1}{2} & \text{if }w\text{ is odd}. 
\end{cases} \label{no-pyt} 
\end{equation} 
\end{lem}
\begin{pf}
Suppose that integers $u,$ $v > 0$ satisfy $u^2-v^2 = -w^2.$ 
This is equivalent to 
\[ v^2-u^2 = w^2,\] 
or equivalently, 
\[ (v+u)(v-u) = w^2.\] 
Let $s = v+u$ and $t = v-u.$ 
The pair $(u,v)$ is given by 
\[ (u,v) = \left(\frac{s-t}{2},\frac{s+t}{2}\right),\] 
and corresponds to the pair $(s,t)$ of integers such that 
\[ st = w^2, \quad s > t > 0, \quad\text{and}\quad s \equiv t \equiv w \pmod 2.\] 
The number of such pairs is equal to the half of the subtraction of $1$ from the number of the positive divisors of $w^2/4$ or $w^2,$ since every pair $(s,t)$ of even divisors of $w^2$ such that $st = w^2$ corresponds to the pair $(s/2,t/2)$ of divisors of $w^2/4$ if $w$ is even. 
By the formula for the number of divisors, the desired number is given by \eqref{no-pyt}.
\end{pf}
\begin{pf}[{\bfseries Proof of Theorem \ref{thm-rat}}]
Let $a = x_1/w,$ $b = y_1/w,$ and $c = z/w.$ 
Then \eqref{star} is equivalent to 
\begin{equation} 
(x_1-z)^2(y_1{}^2+w^2) = (y_1-z)^2(x_1{}^2+w^2). \tag*{$(\star )'$} \label{star-rat} 
\end{equation} 
If $w$ is a multiple of $4$ greater than $4,$ the twice of an odd composite number, or an odd composite number, then the value of \eqref{no-pyt} is greater than $1,$ which implies that $u^2+v^2 = -w^2$ has distinct positive integral solutions $(u,v) = (x_1,x_2),$ $(y_1,y_2).$ 
Letting $x_1{}^2+w^2 = x_2{}^2$ and $y_1{}^2+w^2 = y_2{}^2$ in \ref{star-rat}, 
we obtain 
\[ (x_1-z)^2y_2{}^2 = (y_1-z)^2x_2{}^2\] 
and therefore 
\[ (x_1-z)y_2 = \pm (y_1-z)x_2.\] 
Solving for $z,$ we obtain 
\[ z = \frac{x_1y_2+x_2y_1}{y_2+x_2}, \frac{x_1y_2-x_2y_1}{y_2-x_2}.\] 
Substituting them into $c = z/w$ and considering the sign changes, we obtain the desired rational solutions \eqref{sol-rat}.
\end{pf}

\end{document}